\documentclass[3p,times]{elsarticle}

\usepackage{enumitem}
\usepackage{hyperref}
\usepackage{amsmath}
\usepackage{pifont}
\usepackage[utf8]{inputenc}
\makeatletter
\def\LaTeX{\leavevmode L\raise.42ex
    \hbox{\kern-.3em\size{\sf@size}{0pt}\selectfont A}\kern-.15em\TeX}
\makeatother

\newcommand{\BibTeX}{{\rm B\kern-.05em{\sc
          i\kern-.025emb}\kern-.08em\TeX}}

\makeatletter
\def\@currentlabel{2.1}\label{e:dispaa}
\def\@currentlabel{2.21}\label{e:dispau}
\def\@currentlabel{2.22}\label{e:dispav}
\def\@currentlabel{2.23}\label{e:dispaw}
\def\@currentlabel{2.24}\label{e:dispax}
\def\theequation{\thesection.\@arabic\c@equation}
\makeatother

\renewcommand{\theequation}{\arabic{section}.\arabic{equation}}

\newcommand{\R}{\mathbb R}

\newcommand{\N}{\mathbb N}

\newcommand\supp{\operatorname{supp}}
\newcommand\dist{\operatorname{dist}}

\def \D{\Delta}
\def \O{\Omega}
\def \l{\lambda}
\def \a{\alpha}

\def \o{\omega}
\everymath{\displaystyle}

\newtheorem{thm}{Theorem} [section]
\newtheorem{lem}{Lemma} [section]

\newtheorem{cor}{Corollary} [section]

\newtheorem{rem}{Remark}[section]

\renewcommand{\theequation}{\thesection.\arabic{equation}}
\renewcommand{\thesection}{\arabic{section}}
\renewcommand{\theequation}{\thesection.\arabic{equation}}
\let\ssection=\section\renewcommand{\section}{\setcounter{equation}{0}\ssection}

\begin{document}
\begin{frontmatter}
\title{ Interpolation inequality and some applications}
\author[ah]{Abdellaziz Harrabi}
\ead{abdellaziz.harrabi@yahoo.fr}
\begin{center}
\address{
  Department of Mathematics, Collage of science. Northern Border University, Arar, Saudi Arabia\\Institut Sup\'{e}rieur des Math\'ematiques Appliqu\'ees et de l'Informatique, Universit\`{e} de Kairouan,
  Tunisia}
\end{center}
 \begin{abstract}
  We investigate {\bf explicit} universal estimate of finite Morse index solutions to polyharmonic equations. \,Differently  to previous works \cite{BL2, DDF, fa, H1}, propose here a direct proof using a new interpolation inequality and a delicate boot-strap argument under
    large superlinear and subcritical growth conditions to show that the universal constant grows as a power function of the Morse index.\, Also, our interpolation inequality allows us to provide local $L^p$-$W^{2r,p}$ estimate.
\end{abstract}
 \end{frontmatter}
 \section{ Introduction  }
\subsection{ Interpolation inequalities.}
 Let $n,\,r \geq 2$ be two integer numbers and \,$p\geq 2$ a real number. We designate by $\O$ an open subset of $\R^n$ and $B_R$ the ball of radius $R>0$
 centered at the origin. Let $j=(j_1,j_2,...j_n)$ be a multi index, the weak  $j^{th}$ partial derivative and the magnitude of
the $q^{th}$ gradient of $u\in W^{r,p}_{loc}(\O)$ are respectively defined $a.e$ in $\O$ by
\begin{equation} \label{nab}
 D^ju=\displaystyle\frac{\partial^j u}{\partial x_1^{j_1}...\partial x_n^{j_n}},\,1\leq |j|\leq r \mbox{   and   }|\nabla^{q} u|=\left(\sum\limits_{|j|=q} | D^j u|^p\right)^
 {\frac{1}{p}},\,  1\leq q\leq r.
 \end{equation}
  Let\, $\varepsilon \in (0,1)$\, and \,$1\leq  q\leq r-1.$ From an obvious dilation argument, the standard interpolation inequality \cite{AD} implies

\begin{equation}\label{inter}
 R^{p(q-r)} \int_{B_{R}} |\nabla^q v|^p  \leq  \varepsilon  \int_{B_{R}} |\nabla^{r} v|^p + C \varepsilon^{\frac{-q}{r-q}} R^{-pr} \int_{B_{R}}|v|^p
 ,\quad  v\in W^{r,p}(B_R).
 \end{equation}

   where \, $C=C(n,p,r)$ \,is a positive constant. \,According to \eqref{inter}, one can establish the following weighted interpolation inequality (see \cite{ H1, FH,  RW})
 \begin{equation}\label{ziz}
R^{p(q-r)}\Phi^p_q(v) \leq \varepsilon \Phi^p_r(v)+C \varepsilon^{\frac{-q}{r-q}}R^{-pr}\int_{B_R}|v|^p,
\end{equation}
 where $\Phi_q$ is a family of weighted semi-norms defined by
 \begin{equation*}
\Phi_q(v)= \left(\sup\limits_{0< \alpha<1}(1-\alpha)^{q}\int_{B_{\alpha R}} |\nabla^{q} v|^p\right)^\frac{ 1}{p},\; 0\leq q\leq r.
\end{equation*}
 Inequality \eqref{ziz} together with  the following cut-off function \,$\psi=\psi_{\alpha,R} \in C^{r}_c(\R^n), \alpha \in (0,1)$
 $$ \psi(x)=\exp\left(\left(\frac{\frac{|x|}{R}-\alpha}{\frac{|x|}{R}-\alpha'}\right)^{2r+1}\right)\;\mbox{ if }\, \alpha R<|x|<\alpha'R,\,\psi\equiv 1 \mbox{  if  }\,
  |x|\leq \alpha R \mbox{  and }\psi\equiv 0 \mbox{  if  }|x|\geq \alpha' R \mbox{  where  }\alpha'=\frac{ 1+\alpha}{2}.$$
 are quite useful to provide
the energy estimate which is essential to classify  stable at infinity weak solution  of the $p$-polyharmonic equations \cite{FH} (see also \cite{H1, RW} for $p=2$). \, The
reader may consults \cite{AD, GGS, GT} for further
 applications of  \eqref{ziz}.  When $p\geq 2$ we introduce a new interpolation inequality  which will be
more relevant in providing  integral estimates in various contexts. In particular it will be helpful to establish explicit universal
estimate and local \,$L^p$-$W^{2r,p}$-estimate (see Appendix C).  Moreover, our inequality relies on a more general
cut-off function related to two bounded open subset $\o$ and $\o'$  such that
 $\overline{\o} \subset \o' \subset \overline{\o'} \subset \O.$  Precisely, denote $d= {\rm dist}(\o, \O\backslash\o' ),$ we have
\begin{lem}\label{00}
 There exist $ \psi \in C_c^\infty (\o')$ and a positive constant $C$ depending only on \,$(n,p,k,m)$ such that
 \begin{eqnarray} \label{eq1}
\begin{cases}
 0\leq \psi \leq 1 \mbox{ and }\psi\equiv 1 \;\;\mbox{if} \;\; x \in \o;\\
 |\nabla^k \psi (x)|^p \leq C d^{-k p},\; \forall x \in \o' \mbox{ and } k \in \N.
\end{cases}
\end{eqnarray}
Moreover, we have
\begin{eqnarray}\label{grad}
|\nabla^k \psi^{m}|\leq Cd^{-k} \psi^{(m-k)},\; \;\forall x \in \o' \mbox{ and }m> k.
\end{eqnarray}
\end{lem}
As usual, we used the power function $\psi^m,\;m> r$ as a cut-off function (see \cite{CR, fa,  SZ, H1, FH}). Let $(q,k) \in \N^*\times \N^*,\; q+k=r.$, our main first
 result reads as follows.
\begin{lem}\label{t1}
For every $0<\varepsilon <1,$ there exists a positive constant $C=C(n,r,p,m)$  such that for any $u \in W^{r,p}_{loc}(\O),$ we have
\begin{eqnarray} \label{eq2}
\int_{\o'} |\nabla^q u|^p|\nabla^k (\psi^m)|^p\leq Cd^{-pk}\int_{\o'} |\nabla^q u|^p \psi^{p(m-k)} \leq
\varepsilon \int_{\o'} |\nabla^{r} u|^p\psi^{pm}+C\varepsilon^{1-p^{r}} d^{-pr}\int_{\o'}|u|^p \psi^{p(m-r)}.
\end{eqnarray}
Consequently,
\begin{eqnarray}\label{ij}
  \int_{\o'} |\nabla^{r} u|^p\psi^{pm}\leq 2\int_{\o'} |\nabla^{r}( u\psi^{m})|^p+C d^{-pr}\int_{\o'}|u|^p \psi^{p(m-r)},
 \end{eqnarray}
and
\begin{eqnarray} \label{ij'}
\int_{\o'} |\nabla^q u|^p|\nabla^k \psi^m|^p\leq Cd^{-pk}\int_{\o'} |\nabla^q u|^p \psi^{p(m-k)} \leq \varepsilon \int_{\o'} |\nabla^{r}( u\psi^{m})|^p
+C\varepsilon^{1-p^{r}}  d^{-pr}\int_{\o'}|u|^p \psi^{p(m-r)}.
\end{eqnarray}
\end{lem}
\subsection{ Explicit universal estimate.}
  Consider the following polyharmonic problem:
\begin{equation}\label{k}
 (-\Delta)^r u= f(x,u),  \mbox{ in } \O.
 \end{equation}
 Here, $\O$ is a proper domain of $\R^n$, $u\in C^{2r}(\O)$, $f$ and $f'=\frac{\partial f}{\partial s}$ belong to $C(\Omega\times \mathbb{R})$.
 The associated quadratic form of \eqref{k} is defined by
\begin{equation}\label{qf}
Q_u(h)= \displaystyle\int_{\O}|\textit{D}_{r}h|^2 -\displaystyle\int_{\O}
    f'(x,u)h^2, \;  h\in C_c^{r} (\O),
  \end{equation}
where
$\textit{D}_{r}h = \nabla\Delta^{j-1}h,\,|\textit{D}_{r}h|^2=|\nabla\Delta^{j-1}h|^2  \,\mbox{ if\, $r=2j-1$ } \, \mbox{  and  } \,\textit{D}_{r}h =\D^{j}h,
 \,|\textit{D}_{r}h|^2=(\D^{j}h)^2 \mbox{ if \,$r=2j$ }$, $j\in \N^*.$
\medskip
The Morse index of $u$, denoted by  $i(u)$ is defined as the maximal dimension of all subspaces $V$ of $C_c^{r} (\O)$ such that \\$ Q_{u} (h)< 0,\;\forall\, h \in V\backslash \{0\}$.
In previous works \cite{DDF, fa,  H1}, universal estimate has been established from blow-up technique and some available Liouville-type theorems
 classifying finite Morse index solutions (see also the case of positive solutions in \cite{GS, PQS,  RW,  SZ}.\, However, this procedure fails to derive explicit estimate
 and requires a restrictive asymptotic behaviour condition:\\
$(h_{0}):$ There exists $q>1$ such that  $\displaystyle \lim_{|s|\rightarrow \infty}\frac {f'(x,s)}{|s|^{q-1}}= 1,$ uniformly with respect to $x\in \O$.\\
Thanks to Lemma \ref{t1}, we establish  explicit universal estimate under the following large
 superlinear and subcritical growth conditions:\\  There exist $s_0> 0,$ $c_1>1$ and $1<p_1\leq p_2<\frac{n+2r}{n-2r}$
 such that for all $(x,s)\in \Omega\times\R\backslash [-s_0,s_0],$
\begin{itemize}
 \item [$(h_{1})$] (Super-linearity)\, $f'(x,s)s^{2}\geq p_1f(x,s)s;$
\item [$(h_{2})$]  (Subcritical growth)\, $(p_2+1)F(x,s)\geq f(x,s)s,$  where $F(x,s)=\int_{ 0}^{s}f(x,t)dt;$
\item [$(h_{3})$]  $|(\nabla_y F)(x,s)|\leq c_1(F(x,s)+1)$, \,for all\, $(x,s)\in \Omega\times\R;$
\item [$(h_{4})$]  $|f'(x, s)| \leq c_1$,\, for all \,$(x,s)\in \Omega\times [-s_0,s_0]$, \,\,$|f(x,0)| \leq c_1$ \,and\,$\pm f(x,\pm s_0)\geq \frac{ 1}{c_1}$,\,
 for all\, $x\in \Omega$.
\end{itemize}

 When $f(x,s)=f(s)$ the above assumptions are reduced to $(h_{1})$-$(h_{2})$ (with $\pm f(\pm s_0)>0$)\; and obviously are weaker than $(h_0).$\, Let $K\in C^1(\O)$
 \,be a positive
function such that\,$K,\;|\nabla K| \in L^\infty(\O)$, \,and
  $1<p_1<
 p_2<\frac{n+2r}{n-2r}$\, and denote $s_+=\max(s,0),$ $s_-=\max(-s,0)$ \footnote{If $\O$ is an unbounded domain we assume in addition that $K(x)\geq c_0 >0$ for all $x\in \O.$} . \,The
 nonlinearity  $f(x,s)=K(x)(s_+^{p_2}-s_-^{p_1})$ \,satisfies \,$(h_{1})$-$(h_{4})$\,but violates  $(h_0).$  \,  Let \,$ \alpha \in (0,1),$\, $y \in \O$.\,
Denote \,$\delta_y=\dist(y,\partial \O),$ \,$ d_y=\inf(\alpha,\delta_y)$ \footnote{ In the statement of Theorem \ref{univ}, we used \eqref{nab} with $p=1.$}.  We have
\begin{thm}\label{univ}
Assume that $f$ satisfies $(h_{1})$-$(h_{4})$. \,Then, there exist $\alpha_0\in (0,1),\,\gamma_1>0,\,\gamma_2>0$ and a positive constant
 $C=C(\alpha_0,n,r,p_1,p_2,s_0,c_1)$\, independent
of \,$\O$\, such that for any finite Morse index solution \,$u$ \,of \eqref{k} and for every \,$\alpha\in (0,\alpha_0),$\, we have
  \begin{eqnarray}\label{zz1}
  \sum_{ j=0}^{2r-1}d_y^j|\nabla^j u(y)|\leq C(1+i(u))^{\gamma_2}d_y^{-\gamma_1}, \; \forall \;  y\in \O.
  \end{eqnarray}
Precisely, if\, $\frac{p_2+1}{p_2}<\frac{n}{2r}$ \,then \,$\gamma_1=\frac {4r^2(p_1+1)p_2}{(p_1-1)(2r(p_2+1)-n(p_2-1))}$ \,and
\,$\gamma_2=\gamma_1+\frac {2r(p_2+1)}{2r(p_2+1)-n(p_2-1)}.$
\end{thm}
\begin{rem}
 Denote \,$\O_\alpha=\{ y \in \O, \delta_y \geq \alpha\}$,\,$\alpha\in (0,\alpha_0).$\, As a direct consequence of \eqref{zz1}, we have
 $$\|u\|_{C^{2r-1}(\O_\alpha)}\leq C\alpha^{1-2r-\gamma_1}(1+i(u))^{\gamma_2}\mbox{ and if } y \in \O\backslash \O_\alpha, \mbox{ then }
 \sum_{ j=0}^{2r-1}|\nabla^j u(y)|\leq C(1+i(u))^{\gamma_2}\delta_y^{1-2r-\gamma_1}.$$
\end{rem}
To prove Theorem \ref{univ} we make use of Lemmas \ref{00} and \ref{t1} to obtain a first integral estimate on a ring around $y$
 (see \eqref{pip} in Section 3).\, By virtue of
a variant of the Pohozaev identity \cite{P}, we extend this estimate to a bull centered at $y$ as follows
 \begin{align}\label{hdds}
d_y^{-n}\int_{B(y,\frac {d_y}2)}|f(x,u)|^{\frac{ p_2+1}{p_2}}   \leq C \left(\frac{1+i(u)}{d_y}\right)^{\frac{2(p_1+1)r}{p_1-1}+1}, \;\forall \;y\in \O.
 \end{align}
 As $p_2$ is subcritical, we used a delicate boot strap argument to end the proof of Theorem \ref{univ}. \,Note that estimate \eqref{hdds} holds when
  \,$\frac{n+2r}{n-2r}<p_1\leq p_2,$\,
 but it is not clear which procedure would be helpful to derive \eqref{zz1}\footnote{ Note that the boot strap argument requires a subcritical growth.}. \, Also, inequality
 \eqref{hdds} could be extended to solutions of the $p$-polyharmonic equation .\,However, we do not dispose to any $L^q$-regularity result to star
 the boot strap procedure. \,Regarding the case of bounded domain,
  explicit \,$L^{\infty}$-bounds of finite Morse index solutions of the second order Dirichlet boundary-value problem has been obtained in \cite{HHM,  KA, lec} under similar assumptions
 of \,$(h_1)$-$(h_4)$ which improve the a priori $L^\infty$-estimates stated in \cite{BL2, HRS2}.\,  Also, in  \cite{HMY}  the authors examined the influence of the type boundary conditions
 involving the biharmonic and triharmonic problems to provide similar  explicit \,$L^{\infty}$-bounds. The general higher order case $r\geq 4,$ is more
 difficult since some needed local interior estimates near the boundary are so hard to achieve.\\
This paper is organized as follows: Section 2 is devoted to the proofs of Lemmas \ref{00} and \ref{t1}. In section 3, we give the proof of  Theorem \ref{univ}.  In appendix C,
 we provide the proof of local
 $L^p$-$W^{2r,p}$\, estimate. \\In the following, $C $ (respectively $C_\varepsilon$) denotes always generic positive constants depending only on $(n,p,r,k,m)$ (respectively on
  $(\varepsilon, n,p, r,m)$)
 which could be changed from one line to another.
\section{ Proofs of Lemmas  \ref{00} and \ref{t1}.}
{ \bf Proof of Lemma  \ref{00}.} Set \,$\o_d=\{x \in \O,\; {\rm dist}(x,\o) < \frac d 4\},$\, where \,$d= {\rm dist}(\o, \O\backslash\o' ),$ we have\, $\o \subset \o_d \subset \o'$.
\,Let\, $h =\chi_{\o_d}$ \,be the indicator function of \,$\o_d$ and $g \in C_c^\infty(\R^n)$ a nonnegative function such that
 $\supp(g) \subset B_1 \mbox{ and }\displaystyle\int_{\R^n} g(x) dx = 1.$ Set
 $$g_d(x) = \left(\frac 8{d}\right)^{n}g \left(\frac {8x}{d}\right) \mbox{  and  }
 \psi (x) =\displaystyle\int_{\R^n } g_d(y) h(x-y)dy=\displaystyle\int_{B_{\frac d{8} }} g_d(y) h(x-y)dy.$$
We have $0\leq \psi \leq 1$ and $ \supp(\psi ) \subset \o_d+B_{\frac {d}{8}} \subset \o'$ (see proposition 4.18 in \cite{Br}). Since $\o +B_{\frac {d}{8}} \subset \o_d,$\, then
$ \psi(x)=1 \mbox{ if }x \in \o .$ Also, $\psi \in C_c^\infty(\R^n) \mbox{ with } D^j\psi(x)=\displaystyle\int_{B_{ \frac d{8}}}D^jg_d(y) h(x-y)dy$ (see  proposition 4.20 in
\cite{Br}). Therefore,
  $$|D^j\psi(x)| \leq \displaystyle\int_{B_\l} |D^j g_d|dy \leq (\frac 8{d})^{-|j|}\displaystyle\int_{B_{1}}|D^j g(y)|dy \leq Cd^{-|j|}.$$
   Now, from \eqref{nab} one can see that $ |\nabla^k \psi (x)|^p \leq C d^{-k p},\; \forall x \in \o'\backslash\o \mbox{ and } k \in \N,$ where is $C=C(n,k,p)>0.$\\
  { \bf Proof of   \eqref{grad}.} The proof will be done by working inductively with respect $k\geq 1$. Observe that \eqref{grad} is an immediate consequence of
  \eqref{eq1} if  $k=1$.  Assume now that the following inequality holds for all $1\leq l\leq k$ and $m>l$
 \begin{eqnarray}\label{ass}
|\nabla^l \psi^{m}|\leq Cd^{-l} \psi^{(m-l)},\; \;\forall x \in \o'.
\end{eqnarray}
Let $m> k+1$, fix $j=(j_1, j_2,...,j_n)$ such that $2 \leq |j| \leq k+1$ and  $i_0 \in \{1, \;2, \;...,n\}$ such that $j_{i_0}\neq 0$ and
 denote $j_{-}=(j_1,.., j_{i_0}-1, .. j_n)$. According to Leibnitz's formula, we have
$$\mathcal{D}^{j} \psi^{m}= m \mathcal{D}^{j_{-}}(\psi^{m-1}\frac{\partial \psi }{\partial x_{i_0}})=m\psi^{m-1} \mathcal{D}^{j}
 \psi +m\sum\limits_{\substack{ s+t= j_0\\ t\neq j_0}} a_{s,t}\mathcal{D}^{s}\psi^{m-1}\mathcal{D}^{t}\frac{\partial \psi }{\partial x_{i_0}}, \mbox{ where } |s|+|t|=k, a_{s,t} \in \R. $$  From \eqref{eq1}, we derive
$$|\nabla^k \psi^{m}|\leq C \left(d^{-k}\psi^{m-l}+ \sum_{ 1\leq l\leq k} d^{l-k-1}|\nabla^l \psi^{m-1}|\right), \; \forall x\in \o'.$$
  According to our assumption \eqref{ass}, $m-1>k$ and the above inequality, we derive that \eqref{grad} holds for $k+1.$  This achieves the proof of Lemma \ref{00}.\qed
\subsection{  Proof of Lemma  \ref{t1}.}
 We will use the following elementary inequalities. For $p\geq 2$, $\varepsilon \in (0,1)$, $a,b$ and $c$ positive real numbers , we have
 \begin{eqnarray} \label{y1}
 b^p\leq  2a^{p}+C|a-b|^{p},\,\, \,\,ab^{p-2}c \leq \frac {1}{ p}\varepsilon^{1-p} a^p +\frac {p-2}{ p}\varepsilon b^p +\frac {1}{ p}\varepsilon c^p.
     \end{eqnarray}

 Let $\psi$ the cut-off function defined in Lemma \ref{00} and $m>r$. Inequality \eqref{ij'}  is an immediate consequence of \eqref{eq2} and \eqref{ij}. Also, inequality \eqref{ij} follows
 from \eqref{eq2}.  In fact, from \eqref{nab}, we have
 $|\nabla^r u|^p\psi^{pm} =\displaystyle\sum\limits_{|j|=r} |D^j u|^p\psi^{pm}.$ Thus, the first inequality of \eqref{y1} (with  $a=|D^j ( u\psi^{m})|$ \, and \, $b=|D^j ( u)\psi^{m}|$) and
  Leibnitz's formula \cite{AD} imply
\begin{eqnarray*}
 |\nabla^{r} u|^p\psi^{pm} &\leq &   2|\nabla^{r}( u\psi^{m})|^p+C\displaystyle\sum_{|j|=r}|D^j ( u\psi^{m})-D^j ( u)\psi^{m}|^p \nonumber\\ &\leq &
   2|\nabla^{r}( u\psi^{m})|^p+C\displaystyle\sum_{q+k=r,\,q\neq r}|\nabla^{q} u|^p|\nabla^{k}\psi^{m}|^p.
\end{eqnarray*}
In view of  \eqref{grad}, we get $\displaystyle\int_{\o'}|\nabla^{r} u|^p\psi^{pm}  \leq 2 \displaystyle\int_{\O}|\nabla^{r}( u\psi^{m})|^p+C \displaystyle
 \sum_{q+k=r,\,q\neq r}d^{-pk}\int_{\o'}|\nabla^{q} u|^p\psi^{p(m-k)}.$  Hence,  inequality \eqref{ij} follows from  \eqref{eq2}.

{\bf Proof of \eqref{eq2}.}  Set $I_{q}=d^{-pk}\displaystyle\int_{\o'} |\nabla^q u|_p^p\psi^{p(m-k)}.$ From  \eqref{grad}, we have
 $\displaystyle\int_{\o'} |\nabla^q u|^p|\nabla^k (\psi^m)|^p \leq  C I_q$. Thus, to provide  \eqref{eq2}, we have only to prove the following inequality:
\begin{eqnarray} \label{eq41}
 I_{q} \leq  \varepsilon I_{r}+C\varepsilon^{1-p^{r}} I_{0},\; \forall \;1\leq q\leq r-1.
\end{eqnarray}
We divide the proof of \eqref{eq41} into two steps.

{\bf Step $1$.} We establish the following first-order interpolation inequality:
 \begin{eqnarray} \label{eq4}
 I_{q} \leq  \varepsilon I_{q+1}+C  \varepsilon^{1-p} I_{q-1},\; 1\leq q\leq r-1.
\end{eqnarray}

Recall that $\psi \in C_c^\infty(\o')$ and denote \, $u_{|\o'}$\, the restriction of\, $u$ on\, $\o'$. Observe that by virtue of Meyers-Serrin's
 density theorem \cite{AD} and using Lebesgue's dominated convergence theorem \cite{Br}, one can reduce the proof of \eqref{eq4} to $u_{|\o'}$ belonging to $\in C^r(\o')\cap W^{r,p}(\o').$ \, Let $j=(j_1, j_2, ..., j_n)$ be a multi
 index with $|j|= q\leq r-1$ and
$i_0 \in \{1, \;2, \;...,n\}$ such that $j_{i_0}\neq 0.$  Set $j_{-}=(j_1,.., j_{i_0}-1, .. j_n)$,\, $\;|j_-|=q-1$ and $j_{+}=(j_1,.., j_{i_0}+1, .. j_n),$ \,
$|j_+|=q+1$.  As $p\geq 2$ and $|j|\leq r-1,$ we have
\begin{eqnarray} \label{eq'44}
 |D^j u|^{p-2}D^j u \in C^1(\o')\; \mbox{  and  }\; \frac{\partial (|D^j u|^{p-2}D^j u) }{\partial x_{i_0}}= (p-1)|D^j u|^{p-2}D^{j_+} u. \end{eqnarray}
From \eqref{grad} on has $|\nabla \psi|\leq C d^{-1}$, then integration by parts yields
\begin{eqnarray}\label{eq'4}
d^{-pk}\int_{\o'}  |D^{j}u|^{p}\psi^{p(m-k)}
&=& \;d^{-pk} \int_{\o'}|D^j u|^{p-2}D^j u\frac{\partial D^{j_-}u }{\partial x_{i_0}}\psi^{p(m-k)}
\nonumber\\ &=& \;-(p-1)d^{-pk}\int_{\o'}  |D^{j}u|^{p-2}D^{j_+}uD^{j_-}u\psi^{p(m-k)}\nonumber\\ &
  & -p(m-k)d^{-pk}\int_{\o'} |D^j u|^{p-2}D^j u D^{j_-}u\psi^{p(m-k)-1}\frac{\partial \psi}{\partial x_{i_0}}\nonumber\\ &
\leq & C d^{-pk}\int_{\o'} |\nabla^{q-1}u ||\nabla^{q}u|^{p-2}|\nabla^{q+1}u|\psi^{p(m-k)} \nonumber\\ && +\;C d^{-(pk+1)} \int_{\o'}|\nabla^{q-1}u ||\nabla^{q}u|^{p-1} \psi^{p(m-k)-1}.
\end{eqnarray}
Taking into account that $I_q=\sum\limits_{|j|=q} d^{-pk}\int_{\o'}  |D^{j}u|^{p}\psi^{p(m-k)}$ \,with\, $k=r-q,$ \, so inequality \eqref{eq'4} implies
 \begin{eqnarray} \label{ma}
 I_{q} \leq  C(J_1+J_2) \mbox{ where } J_1=  d^{-pk}\int_{\o'}|\nabla^{q-1} u| |\nabla^q u|^{p-2}|\nabla^{q+1} u|\psi^{p(m-k)} \mbox{ and}\;\; J_2= d^{-(pk+1)} \int_{\o'}|\nabla^{q-1} u| |\nabla^q u|^{p-1}\psi^{p(m-k)-1}.
 \end{eqnarray}

Observe that $pk= (k+1)+k(p-2)+(k-1),$ \,$p(m-k)= (m-(k+1))+(p-2)(m-k)+(m-(k-1))$ \,(respectively $pk+1=(k+1)+(p-2)k+k $ \,and \,$p(m-k)-1=(m-(k+1))+(p-2)(m-k)+(m-k)$. \,Thus, inequality \eqref{y1} with \,$a=d^{-(k+1)}|\nabla^{q-1} u| \psi^{m-(k+1)},\;b=d^{-k}|\nabla^{q} u| \psi^{(r-k)} \mbox{ and } c=d^{-(k-1)}|\nabla^{q+1} u| \psi^{m-(k-1)} (\mbox{ respectively }c=d^{-k}|\nabla^{q} u| \psi^{(r-k)}), $ implies
\begin{eqnarray*}
 J_1\leq \frac{1}{p}\varepsilon^{1-p} I_{q-1}+\frac{(p-2)}{p}\varepsilon I_{q}+\frac{1}{p}\varepsilon  I_{q+1} \mbox{ and }J_2\leq  \frac{1}{p} \varepsilon^{1-p}  I_{q-1}+\frac{p-1}{p}\varepsilon I_{q}.
\end{eqnarray*}
 Combining the above inequalities with \eqref{ma} we deduce $(1-2C\varepsilon) I_{q} \leq C \varepsilon^{1-p}  I_{q-1} +C\varepsilon I_{q+1}.$
Hence, the inequality \eqref{eq4} follows by replacing $ \varepsilon $\, by \,$\frac {\varepsilon}{4(1+C)}.$\\
{\bf Step $2$. End of the proof of \eqref{eq41}.} The case $r=2$, or $r\geq 3$ and $q=1$ are an immediate consequence of \eqref{eq4}.  Let $r\geq 3,$ $2\leq q\leq r-1$ and
$2\leq t\leq q$ and set $S_t=\sum_{i=2}^{t} I_{i}.$ We apply \eqref{eq4} where
one substitutes  $q$ by $t-i$ and $\varepsilon$ by $\varepsilon^{p^i}$, we derive $C^i\varepsilon^{-p^i} I_{t-i}\leq C^{i+1} \varepsilon^{-p^{i+1}} I_{t-i-1} +  C^i I_{t-i+1}.$  Since $S_{t}\leq S_q$ and  $0<\varepsilon <1$, the summation of the above inequalities from \, $i=0$\, to\, $i=t-1$\, yields

\begin{eqnarray}\label{1.33}
I_{t} \leq  C \varepsilon^{1-p^{r}} I_{0} + \varepsilon I_{t+1}+C \varepsilon S_{q}\; \mbox{  if   } 2\leq t\leq q.
 \end{eqnarray}
Summing now \eqref{1.33} from \,$t=2$ to $t=q$ and substituting  $\varepsilon$ by $\frac{ \varepsilon}{2(C+1)},$ we arrive at $S_q\leq C\varepsilon^{1-p^{r}} I_{0} +
 \varepsilon I_{q+1}, \mbox{ for all } 1\leq q\leq r-1.$   Combining  \eqref{1.33} with $t=q$ and the last inequality, we obtain
 \begin{eqnarray}\label{dawn}
I_{q}\leq C \varepsilon^{1-p^{r}} I_{0}+\varepsilon I_{q+1},\; 1 \leq  q\leq r-1.
\end{eqnarray}
 To end the proof of \eqref{eq41}, we iterate \eqref{dawn} as follows
 \begin{eqnarray}\label{1.3}
\left\{\begin{array}{llllll}
I_{q} \leq  C \varepsilon^{1-p^{r}} I_{0} +  I_{q+1},\\
  I_{q+1} \leq +C \varepsilon^{1-p^{r}} I_{0,r} + I_{q+2},\\
\vdots \\
 I_{r-1} \leq C \varepsilon^{1-p^{r}} I_{0} + \varepsilon I_{r}.\\ \end{array} \right.
 \end{eqnarray}
Hence, the summation of the above inequalities yields  $I_{q}\leq C \varepsilon^{1-p^{r}} I_{0}+\varepsilon I_{r},$
which is the desired inequality \eqref{eq41}. The proof of Lemma \ref{t1} is completed.\qed
 \section{ Proof of Theorem \ref{univ}.}
\subsection{Preliminary  results.}
  $B(y,\l)$ stands for the ball of radius $\l>0$ centered at $y\in \R^n$.  \, Let $\psi$ be the cut-off function defined in Lemma \ref{00} related to two
   open subset \,$\o$ \,and \,$\o'$ of $B(y,\l).$
    Thanks to Lemma \ref{t1} with $p=2$, we establish the following technical lemma:
\begin{lem}\label{L2}
For every $0<\varepsilon<1,$ there exists a positive constant $ C_\varepsilon=C(n,m,r,\varepsilon)$ such that, for all $u\in H^{r}(B(y,\l)),$ we have
 \begin{align}
\label{T488}
  \int_{B(y,\l)}\big|| \textit{D}_{r} (u\psi^{m})|^2- \textit{D}_{r}u  \textit{D}_{r}(u\psi^{2m})\big|
\leq  \varepsilon \int_{B(y,\l)}|\nabla^{r} (u\psi^{m})|^{2}+C_\varepsilon  d^{-2r}\int_{B(y,\l)} |u|^2\psi^{2(m-r)};
\end{align}
\begin{align}
\label{T477}
  \int_{B(y,\l)}\big|| \textit{D}_{r} u|^{2} \psi^{2m}-\textit{D}_{r}u  \textit{D}_{r}(u\psi^{2m})\big|\leq
   \varepsilon \int_{B(y,\l)}|\nabla^{r} (u\psi^{m})|^{2}+C_\varepsilon d^{-2r}\int_{B(y,\l)} |u|^2\psi^{2(m-r)};
\end{align}

\begin{align}\label{T99}
 \int_{B(y,\l)} | \nabla^{r}(u\psi^{m}) |^2
  \leq C\left(\int_{B(y,\l)}| \textit{D}_{r} u|^{2} \psi^{2m}+d^{-2r}\int_{B(y,\l)} |u|^2\psi^{2(m-r)}\right);
\end{align}
where $d= {\rm dist}(\o, B(y,\l)\backslash\o' ).$ Moreover, if $u\in H^{r+1}(B(y,\l)),$ we have
\begin{align}\label{T478}
 \int_{B(y,\l)}\big|
 \textit{D}_{r}u\cdot\textit{D}_{r}(\nabla u\cdot (x-y))\psi^{2m}-\textit{D}_{r}u \cdot\textit{D}_{r}(\nabla u\cdot (x-y)\psi^{2m})\big|\leq
   C(1+(\frac{\l }{d})^2)\left( \sum_{1\leq q\leq r} d^{-2(r-q)}\int_{\o'\backslash\o} |\nabla^{q}u|^2+ d^{-2r}\int_{\o'\backslash\o} |u|^2\right).
\end{align}
\end{lem}
 {\bf Proofs of Lemma \ref{L2}.}\\
 For $v\in H^{r}(B(y,\l))$ and $\eta \in C_c^r(B(y,\l))$, set $A(\eta,v):=\textit{D}_{r}( v\eta)-\eta \textit{D}_{r} v$ \footnote{  Both $ \textit{D}_{r}$ and $A(\eta,v)$
  are respectively scalar operators if $r$ is even, and $n-$vectorial operators if $r$ is odd.}.  A simple computations yield
  $$|\textit{D}_{r} (u\eta)|^2 -\eta^{2}| \textit{D}_{r} u|^{2}=2\eta \textit{D}_{r} u\cdot A(\eta,u)+ |A(\eta,u)|^2
,\,\,\eta^{2}\textit{D}_{r} u \cdot \textit{D}_{r}v-
  \textit{D}_{r}u \cdot \textit{D}_{r}(v\eta^{2})=-  \textit{D}_{r} u\cdot A(\eta^2,v)
,$$ and $|A(\eta,v)|\leq C\sum_{q+k=r,\,
q\neq r} |\nabla^{q}v||\nabla^{k}\eta|.$ Therefore,
\begin{align*}
\big| | \textit{D}_{r} (u\eta)|^2- \eta^{2}| \textit{D}_{r} u|^{2}\big|
 \leq C \sum_{q+k=r,\, q\neq r} \left(| \textit{D}_{r} u||\nabla^{q}u|( \eta|\nabla^{k}\eta|+|\nabla^{k}(\eta^2)|)
 +|\nabla^{q}u|^2 |\nabla^{k}\eta|^2\right).
\end{align*}
  and
 \begin{align*}
| \big|\eta^{2}\textit{D}_{r} u \cdot \textit{D}_{r}v-
  \textit{D}_{r}u \cdot \textit{D}_{r}(v\eta^{2})\big|
  \leq C| \textit{D}_{r} u|\sum_{q+k=r,\, q\neq r} |\nabla^{q}v||\nabla^{k}(\eta^2)|.
\end{align*}
Choosing now $\eta=\psi^{m}$ and using \eqref{grad}, we obtain\footnote{ Observe that\, $ |\textit{D}_{r} u|\leq C|\nabla^{r}u|$.}
\begin{align}
\label{qu2}
\int_{B(y,\l)}\big|\left( \psi^{2m}\textit{D}_{r} u \cdot \textit{D}_{r}v-
  \textit{D}_{r}u \cdot \textit{D}_{r}(v \psi^{2m}) \right)\big|
\leq  C S_1(u,v);
\end{align}
\begin{align}
\label{qu1}
\int_{B(y,\l)}\big|| \textit{D}_{r} (u\psi^{m})|^2-\psi^{2m}|\textit{D}_{r}|)\big| \leq  C (S_1(u,u)+S_2(u));
\end{align}

$$ \mbox{ where }S_1(u,v)=\int_{\o'\backslash\o}|\nabla^r u|\left(\sum_{0\leq q\leq r-1} d^{r-q}|\nabla^{q}v|\psi^{2m+q-r}\right)\mbox{ and }
S_2(u)= \int_{\o'\backslash\o}\left(\sum_{0 \leq q \leq r-1} d^{-2(r-q)}|\nabla^{q}u|^2\psi^{2m-2(r-q))}\right).$$
We invoke Cauchy-Schwarz's inequality:
$$\left|a_r\sum_{q=0}^{r-1}a_q \right| \leq \epsilon |a_r|^2+C_{\epsilon} \sum_{q=0}^{r-1}|a_q|^{2},\,(a_0, a_1,...   a_r) \in \R^{r+1},$$
with $a_r=|\nabla^r u|\psi^{m}$ \,and $a_q=d^{q-r}|\nabla^q u|\psi^{(m+q-r)}$ if $q=0,1...r-1$.  As  $pm-(r-q)=(p-1)m+(m-(r-q))$, we arrive at
\begin{align}
\label{quv}
S_1(u,u)+S_2(u)  \leq  +\varepsilon \int_{B(y,\l)}|\nabla^{r} u|^{2}\psi^{2m}+C_\varepsilon  \sum_{q+k=r, q\neq r} d^{-2k}\int_{B(y,\l)} |\nabla^{q}u|^2|\psi^{2(m-k)}.
 \end{align}
  Hence, inequality \eqref{T488} follows from  \eqref{qu2} (with $v=u$), \eqref{quv} and inequalities \eqref{ij}, \eqref{ij'} of Lemma \ref{t1}.\\
Collecting now , inequalities \eqref{qu1}, \eqref{quv}, \eqref{ij} and \eqref{ij'}, we obtain
\begin{align}
\label{T476}
  \int_{B(y,\l)}\big|| \textit{D}_{r} u|^{2} \psi^{2m}| \textit{D}_{r} (u\psi^{m})|^2-\big|\leq
   \varepsilon \int_{B(y,\l)}|\nabla^{r} (u\psi^{m})|^{2}+C_\varepsilon d^{-2r}\int_{B(y,\l)} |u|^2\psi^{2(m-r)}.
\end{align}
So, inequality \eqref{T477} is an immediate consequence of \eqref{T488} and \eqref{T476}.\\
{\bf Proof of \eqref{T99}.} We appeal to the following higher order Calderon-Zygmund's inequality ( see the proof in the appendix of \cite{CP}, see also \cite{GGS}):
\begin{equation}\label{cazyg}
\int_{B(y,\l)} |\nabla^r v|^2 \leq C \int_{B(y,\l)} |\textit{D}_{r} v|^2 , \forall\; v\in H_0^{r}((B(y,\l)),
\end{equation}
where $C$ is a positive constant depending only on $(n,r,p).$ Therefore, combining \eqref{cazyg} (with $v=u\psi^{m}$) and \eqref{T476}, we provide \eqref{T99}.\\
{\bf Proof of \eqref{T478}.} As $$S_1(u,v)=\int_{\o'\backslash\o}\left(|\nabla^r u|\sum_{0\leq q\leq r-1} d^{r-q}|\nabla^{q}v|\psi^{2m+q-r}\right),$$
  the above Cauchy-Schwarz's inequality, yields
\begin{align}
S_1(u,v)  \leq  \int_{\o'\backslash\o}|\nabla^{r} u|^{2}+ C\sum_{0\leq q\leq r-1} d^{-2(r-q)}\int_{\o'\backslash\o} |\nabla^{q}v|^2.
\end{align}
Fix now $v =\nabla u \cdot(x-y)$ and taking into account that\\
$ |\nabla^q (\nabla u \cdot(x-y))|^{2}\leq C(\l^2  |\nabla^{q+1} u |^{2}+ |\nabla^{q} u |^{2})\leq C(\frac{\l }{d})^2d^2  |\nabla^{q+1}  u |^{2}+ |\nabla^{q}  u |^{2})$,
we deduce that
\begin{align}
\label{quvv}
S_1(u,v)  \leq C(1+(\frac{\l }{d})^2)\left( \sum_{1\leq q\leq r} d^{-2(r-q)}\int_{\o'\backslash\o} |\nabla^{q}u|^2\right)+ d^{-2r}\int_{\o'\backslash\o} |u|^2.
\end{align}
Thus, the  proof of \eqref{T478} follows by collecting \eqref{qu2}, \eqref{quvv}. This ends the proofs of Lemma \ref{L2}.\qed
At last, in view of assumptions  $(h_{1})$-$(h_{4}),$ we have (see the proof in Appendix A):

 \begin{lem}\label{L11} Let  $t> 1$ and set $q_1=\frac {p_2+1}{ p_2}.$ There exists a positive constant $C=C(s_0,p_1, p_2,c_1)$ such that for all $(x,s)\in \O\times\R$, we have
 \begin{itemize}
 \item  $[1]$   $f'(x,s)s^2\geq p_1f(x,s) -C;$
 \item  $[2]$   $ (p_2+1)F(x,s)\geq f(x,s)s -C$;
 \item  $[3]$   $ |s|^{p_1+1} \leq C(|f(x,s)s|+1)$,  \,  $|f(x,s)s|\leq  f(x,s)s+C$ and $|F(x,s)|\leq  C(f(x,s)s+1);$
 \item  $[4]$     $ |f(x,s)|^{q_1}  \leq C(|f(x,s)s| +1)$ and $ |f(x,s)|^{\frac {t}{ p_2}}  \leq C(|s|^t +1);$
 \item  $[5]$   For all $\varepsilon \in\, (0,1)$, $0\leq a\leq 1$ and $b>0$ we have $as^2b \leq C+\varepsilon |f(x,s)s|a^{\frac{p_1+1}2}+
  \varepsilon^{\frac{-2}{p_1-1}} b^{\frac{p_1+1}{p_1-1}}.$
 \end{itemize}
 \end{lem}

 \subsection{ End of the proof of Theorem \ref{univ}.}
Recall that $ d_y=\inf(\alpha,\delta_y),$ where $\delta_y=\dist(y,\partial \O)$, $y \in \O$  and $\alpha \in (0,1).$
 For $j=1, 2,\cdots, i(u)+1,$  set $$ A_j:= \{x\in \mathbb{R}^{n};\; a_j<|x-y|<b_j\},\,a_j =\dfrac{2(j+i(u))}{4(i(u)+1)}d_y;
\;\;b_j =\dfrac{2(j+i(u))+1}{4(i(u)+1)}d_y  \mbox{ and } $$
$$ A'_j:= \{x\in \mathbb{R}^{n};\; a'_j<|x-y|<b'_j\},\;a'_j =\dfrac{2(j+i(u))-\frac{1}{2}}{4(i(u)+1)}d_y.
\;\;b'_j =\dfrac{2(j+i(u))+\frac{3}{2}}{4(i(u)+1)}d_y.$$
Observe that $\subset \overline {A_j}\subset A'_j\subset \overline {A_j} \subset B(y,d_y)$ and let $ \psi_j \in C_c^r(B(y,d_y)$ be the cut-off function defined in Lemma\ref{00}
with $\o=A_j$ and $\o'=A'_j$ and satisfying  $\supp(\psi_{j})
\subset A'_j,\,0\leq \psi_{j}\leq 1 \mbox{ if } x\in A'_j \mbox{ and } \psi_{j}=1 \mbox{ if } x\in A_j.$ Moreover, we have
\begin{eqnarray}\label{gradj}
\begin{cases}
 {\rm dist}(A_j, B(y,d_y)\backslash A'_j )=\dfrac{d_y}{2(i(u)+1)} \mbox{ and  \eqref{grad} implies },\\
 |\nabla^k(\psi^m_{j})(x)|^2\leq C\psi_{j}^{2(m-k)} \left(\frac{1+i(u)}{d_y}\right)^{2k},\; \forall x \in B(y,d_y).
\end{cases}
\end{eqnarray}
 From inequality \eqref{ij'} we derive
\begin{equation}\label{qnf}
\sum_{1\leq q\leq r-1} d^{-p(r-q)}\int_{A_{j}} |\nabla^{q}u|^p|\leq \varepsilon \int_{B(y,d_y)}|\nabla^{r} (u\psi_{j}^{m})|^{2}+C_\varepsilon \left(\frac{1+i(u)}{d_y}\right)^{2r}
\int_{B(y,d_y)} u^2\psi_{j}^{2(m-r)}.
\end{equation}
In the sequel we choose  $m=\frac{(p_1+1)r }{2}$ so that $m> r$ and $\frac{(p_1+1)(m-r)}{p_1-1}=m.$ Thus, point 5 of Lemma \ref{L11} with $s= u,\; a=\psi_{j_0}^{2(m-r)}\mbox{ and }
b=(\frac{1+i(u)}{d_y})^{-2r},$ yields
\begin{eqnarray}\label{uf}
 \left(\frac{1+i(u)}{d_y}\right)^{-2r} \int_{B(y,d_y)} u^2\psi_{j_0}^{2(m-r)} \leq Cd_y^n+ \varepsilon\int_{B(y,d_y)} |f(x,u)u|\psi_{j_0}^{2m}+
 C_ \varepsilon d_y^n\left(\frac{1+i(u)}{d_y}\right)^{\frac{2(p_1+1)r}{p_1-1}}.
\end{eqnarray}
Next observe that $supp (u\psi_{j}^{m}) \cap supp (u\psi_{l}^{m})= \emptyset,\forall \;1 \leq l\neq j\leq 1+i(u),$ then according to the definition of the
   quadratic form  \eqref{qf} we derive
  $$Q_u\left(\sum_{1}^{1+i(u)} \l_ju\psi_{j}^{m}\right)=\sum_{1}^{1+i(u)}\l_j^2Q_u(u\psi_{j}^{m}).$$ in view of  the definition of $i(u)$, there exists
  $j_{0}\in \{1,2,...,1+i(u)\}$ such that $Q_u(u\psi_{j_0}^{m})\geq 0.$ Therefore, point 1 of
   Lemma \ref{L11} implies
\begin{equation}\label{qf}
p_1\displaystyle\int_{B(y,d_y)}f(x,u)u\psi_{j_0}^{2m} -Cd_y^n \leq \displaystyle\int_{B(y,d_y)}f'(x,u)u^2\psi_{j_0}^{2m}\leq  \displaystyle\int_{B(y,d_y)}|\textit{D}_{r} (u\psi_{j_0}^{m})|^2.
  \end{equation}
We divide the proof into three steps.\\
{\bf Step 1.} We shall prove the following estimate
 \begin{align}\label{pip}
\begin{split}
&\sum_{1\leq q\leq r} d^{-2(r-q)}\int_{A_{j_0}} |\nabla^{q}u|^2 +\int_{A_{j_0}} |f(x,u)u|
\leq  C d_y^n\left(\frac{1+i(u)}{d_y}\right)^{\frac{2(p_1+1)r}{p_1-1}}.
\end{split}
\end{align}
Multiplying equation \eqref{k} by $-\frac{1+p_1}2u \psi_{j_0}^{2m},$ integrating by parts, we obtain
$$-\frac{1+p_1}2\displaystyle\int_{B(y,d_y)}f(x,u)u\psi_{j_0}^{2m}=-\frac{1+p_1}2\displaystyle\int_{B(y,d_y)}\textit{D}_{r}u  \textit{D}_{r}(u\psi_{j_0}^{2m}).$$
We combine the last equality with \eqref{qf} and point 3 of  Lemma \ref{L11}, yields
\begin{equation*}
\frac{p_1-1}2\left(\displaystyle\int_{B(y,d_y)} |f(x,u)u|\psi_{j_0}^{2m}+ \displaystyle\int_{B(y,d_y)}|\textit{D}_{r} (u\psi_{j_0}^{m})|^2\right)
\leq Cd_y^n +\frac{p_1+1}2 \displaystyle\int_{B(y,d_y)}\left(|\textit{D}_{r} (u\psi_{j_0}^{m})|^2 -\textit{D}_{r}u \textit{D}_{r} (u\psi_{j_0}^{2m})\right).
  \end{equation*}
It follows from  \eqref{T488} and \eqref{cazyg} that
\begin{equation*}
\int_{B(y,d_y)}|\nabla^{r}( u\psi_{j_0}^{m})|^2 +\int_{B(y,d_y)} |f(x,u)u|\psi_{j_0}^{2m}
\leq  Cd_y^n +\varepsilon \int_{B(y,d_y)}|\nabla^{r} (u\psi_{j_0}^{m})|^{2}+C_\varepsilon \left(\frac{1+i(u)}{d_y}\right)^{2r} \int_{B(y,d_y} u^2\psi_{j_0}^{2(m-r)}.
\end{equation*}
Collecting the last inequalities with \eqref{qnf}  and \eqref{uf}, we get \footnote{ Observe that \,$d_y^n\leq d_y^n\left(\frac{1+i(u)}
{d_y}\right)^{\frac{2(p_1+1)r}{p_1-1}}$ as $d_y=\inf(\alpha, \delta_y) <1$.}
\begin{align*}
\begin{split}
&\sum_{1\leq q\leq r} d^{-p(r-q)}\int_{A_{j_0}} |\nabla^{q}u|^p| +\int_{B(y,d_y)} |f(x,u)u|\psi_{j_0}^{2m}
\leq  C d_y^n\left(\frac{1+i(u)}{d_y}\right)^{\frac{2(p_1+1)r}{p_1-1}}.
\end{split}
\end{align*}
Therefore, inequality \eqref{pip} follows as $\psi_{j_0}(x)=1$\, if \,$x\in A_{j_0}.$\\
{\bf Step 2.}  We shall use the following identity (see the proof in appendix B):
\begin{align}\label{bbcaa}
  \textit{D}_{r}u \textit{D}_{r}(\nabla u \cdot (x-y))=\frac{1}{2}\nabla(|\textit{D}_{r} u|^2)\cdot (x-y)+r|\textit{D}_{r} u|^2,
\end{align}
to establish a variant
 of the Pohozaev identity  and we exploit \eqref{pip} to prove that:
\begin{align}\label{dds}
d_y^{-n}\int_{B(y,\frac {d_y}2)}|f(x,u)|^{q_1}   \leq C \left(\frac{1+i(u)}{d_y}\right)^{\frac{2(p_1+1)r}{p_1-1}+1}.
 \end{align}
 Recall that $ A_{j_0}= \{x\in \mathbb{R}^{n};\; a_{j_0}<|x-y|<b_{j_0}\}$.  We invoke again Lemma \ref{00} with $\o=B(y, a_{j_0})$,\, $\o'=B(y, b_{j_0})$ and let
 $\psi\in C_c^r( B(y,b_{j_0}))$ such that
\begin{eqnarray}\label{mz}
\psi\equiv 1  \mbox{ for all } x\in B(y, a_{j_0}) \mbox{ and } |\nabla^k
\psi^{2m}|\leq  C\left(\frac{1+i(u)}{d_y}\right)^k \psi^{2m-k} \;\forall x\in A_{j_0} \mbox{ and } k=1,2,..,r.
\end{eqnarray}
 Multiplying equation \eqref{k} by $u \psi^{2m}$ (respectively by  $(\nabla u\cdot (x-y)) \psi^{2m}$) and integrating by parts, we get
  \begin{eqnarray*}
\displaystyle\int_{B(y,d_y)} \textit{D}_{r}u \textit{D}_{r}(u\psi^{2m}) = \displaystyle\int_{B(y,d_y)}f(x,u)u \psi^{2m} \mbox{ respectively }\displaystyle\int_{B(y,d_y)}
\textit{D}_{r}u \textit{D}_{r}(\nabla u\cdot (x-y) \psi^{2m}) = \displaystyle\int_{B(y,d_y)}f(x,u)\nabla u\cdot (x-y)) \psi^{2m}.
\end{eqnarray*}
According to inequality \eqref{T477} (respectively \eqref{T478}, \eqref{mz} and \eqref{pip}), we derive

 \begin{align}\label{mb}
\int_{B(y,d_y)}| \textit{D}_{r} u|^{2} \psi^{2m}-\int_{B(y,d_y)}f(x,u)u\phi^{2m}\leq
   \varepsilon \int_{B(y,d_y)}|\nabla^{r} (u\psi^{m})|^{2}+C_\varepsilon  (\frac{1+i(u)}{d_y})^{2r}\int_{B(y,d)} u^2\psi^{2(m-r)};
\end{align}

\begin{eqnarray*}
\displaystyle\int_{B(y,d_y)}  \psi^{2m}\textit{D}_{r}u \textit{D}_{r}(\nabla u\cdot (x-y))\leq\displaystyle\int_{B(y,d_y)}f(x,u)\nabla u\cdot (x-y)\psi^{2m}  +
C d_y^n\left(\frac{1+i(u)}{d_y}\right)^{\frac{2(p_1+1)r}{p_1-1}}+ C \left(\frac{1+i(u)}{d_y}\right)^{-2r} \int_{A_{j_0}} u^2\psi^{2(m-r)}.
\end{eqnarray*}
As above, using  point 5 of Lemma \ref{L11} and \eqref{pip}, there holds that
\begin{eqnarray}\label{uf'}
C \left(\frac{1+i(u)}{d_y}\right)^{-2r} \int_{A_{j_0}} u^2\psi_{j_0}^{2(m-r)} \leq Cd_y^n+\frac 12 \int_{A_{j_0}} |f(x,u)u|\psi_{j_0}^{2m}+
 Cd_y^n\left(\frac{1+i(u)}{d_y}\right)^{\frac{2(p_1+1)r}{p_1-1}}\leq   C d_y^n\left(\frac{1+i(u)}{d_y}\right)^{\frac{2(p_1+1)r}{p_1-1}}.
\end{eqnarray}
Combining these inequalities we get
 \begin{eqnarray}\label{nzz}
\displaystyle\int_{B(y,d_y)}  \psi^{2m}\textit{D}_{r}u \textit{D}_{r}(\nabla u\cdot (x-y))\leq\displaystyle\int_{B(y,d_y)}f(x,u)\nabla u\cdot (x-y)\psi^{2m}  +
C d_y^n\left(\frac{1+i(u)}{d_y}\right)^{\frac{2(p_1+1)r}{p_1-1}}.
\end{eqnarray}
In one hand, integration by parts of  the first term of the right hand-side, gives
 \begin{eqnarray*}
\displaystyle\int_{B(y,d_y)}f(x,u)\nabla u\cdot (x-y)\psi^{2m} =-n\int_{B(y,d_y)} F(x,u)\psi^{2m}-\int_{A_{j_0} }
 F(x,u)(\nabla \psi^{2m}\cdot (x-y) )
  +\int_{B(y,d_y)}(\nabla_x F)(x,u)\cdot (x-y)\psi^{2m},
\end{eqnarray*}
 Invoking now assumption $(h_3)$ with points 2-3 of Lemma \ref{L11}, \eqref{mz} and  using again \eqref{pip}, it follows that \footnote{We also use that $|x-y|\leq  d_y \leq \alpha.$}  imply

\begin{eqnarray}\label{mzz}
\displaystyle\int_{B(y,d_y)}f(x,u)(\nabla u\cdot (x-y)\psi^{2m} =(C\a+-\frac{n}{p_2+1})\int_{B(y,d_y)} f(x,u)u\psi^{2m} + Cd_y^n(1+i(u))
 \left(\frac{1+i(u)}{d_y}\right)^{\frac{2(p_1+1)r}{p_1-1}}.
\end{eqnarray}

On the other hand, using \eqref{bbcaa} and integrating by parts we derive
\begin{eqnarray*}
\displaystyle\int_{B(y,d_y)}  \psi^{2m}\textit{D}_{r}u \textit{D}_{r}(\nabla u\cdot (x-y))=\frac{2r-n}{2}\int_{B(y,d_y)}|\textit{D}_{r} u|^2\psi^{2m} -
\frac{1}{2}\int_{A_{j_0}}|\textit{D}_{r}u|^2(\nabla\psi^{2m}\cdot (x-y)).
\end{eqnarray*}
As $|\textit{D}_{r} u|^2\leq |\nabla^{r} u|^2$ and $|x-y|\leq 1$, it follows from \eqref{mz} and \eqref{pip}  that
\begin{eqnarray*}
\displaystyle\int_{B(y,d_y)}  \psi^{2m}\textit{D}_{r}u \textit{D}_{r}(\nabla u\cdot (x-y))=\frac{2r-n}{2}\int_{B(y,d_y)}|\textit{D}_{r} u|^2\psi^{2m} +Cd_y^n
 \left(\frac{1+i(u)}{d_y}\right)^{\frac{2(p_1+1)r}{p_1-1}+1}.
\end{eqnarray*}
 Collecting  inequalities \eqref{nzz}, \eqref{mzz} and the last equality we arrive at
\begin{eqnarray}\label{ pyt}
\left(\frac{2n}{(p_2+1)(n -2r)}-\frac{C\alpha}{(n -2r}\right)\int_{B(y,d_y)} f(x,u)u\psi^{2m}- \int_{ B(y,d_y) }|\textit{D}_{r} u|^2\psi^{2m}
\leq Cd_y^n\left(\frac{1+i(u)}{d_y}\right)^{\frac{2(p_1+1)r}{p_1-1}+1}.
\end{eqnarray}
We choose $\alpha=\alpha_0 \in (0,1)$ small enough so that $\frac{2n}{(p_2+1)(n -2r)}-\frac{C\alpha_0}{(n -2r}>1$ and we combine
 the above inequality with \eqref{mb} and \eqref{T99}, we deduce that \footnote{ Recall that $\frac{2n}{(p_2+1 )(n -2r)}>1.$}
  \begin{eqnarray*}
\int_{B(y,d_y)}|\nabla^{r}u|^2\psi^{2m}+\int_{B(y,d_y)}f(x,u)u\psi^{2m}
&\leq& \varepsilon \int_{B(y,d_y)}|\nabla^{r} (u\psi^{m})|^{2}
+ C_\varepsilon  \frac{(1+i(u))^{2r}}{d^{2r}_y}\int_{B(y,d)} u^2\psi^{2(m-r)}\nonumber\\ && + \; Cd_y^n\left(\frac{1+i(u)}{d_y}\right)^{\frac{2(p_1+1)r}{p_1-1}+1}.
\end{eqnarray*}
Inequality \eqref{uf'} and points 3 of Lemma \ref{L11} imply
\begin{align*}
 \int_{B(y,d_y)}|\nabla^{r}(u\psi^{m})|^2+ \int_{B(y,d_y)}|f(x,u)u|\psi^{2m}\leq
    C  d_y^n\left(\frac{1+i(u)}{d_y}\right)^{\frac{2(p_1+1)r}{p_1-1}+1}.
 \end{align*}
 Observe now that  $\psi\equiv 1$ for all $ x\in B(y,\frac {d_y}2)\subset B(y,a_{j_0}),$  so estimate \eqref{dds} follows from the above inequality and point 4 of Lemma \ref{L11}.

{\bf Step $3$. Boot-strap procedure.} Set $\lambda =\frac {d_y}2 <1$, \, $u_\l(x)=u(y+\lambda x)$ and $g_\l(x)=f(y+\l x,u(y+\l x)),\; x \in B_1,$ then $u_\l$  satisfies
\begin{eqnarray}\label{3.4}
  (-\Delta u_\l)^r=\l^{2r}g_\l \mbox{ in } B_1.
\end{eqnarray}
By virtue of \eqref{dds}, we have
\begin{align}\label{ddss}
\int_{B_1}|g_\l|^{q_1} =2^n d_y^{-n}\int_{B(y,\frac {d_y}2)}|f(x,u)|^{q_1}  \leq C\left(\frac{1+i(u)}{d_y}\right)^{\frac{2(p_1+1)r}{p_1-1}+1}.
 \end{align}
  We invoke local $L^p$-$W^{2r,p}$ estimate (see  Corollary \ref{tloc} in the Appendix C) and Rellich-Kondrachov's theorem \cite{GGS}.
   Let  $q>1$, then point 3 of Lemma \ref{L11} implies \footnote{ Observe that $\lambda =\frac {d_y}2 <1.$}
\begin{eqnarray}\label{3.7}
\|u_\l \|_{L^{q^*}(B_\frac 12)}\leq C\|u_\l\|_{W^{2r,q}(B_\frac 12))} \leq C ( \|g_\l \|_{L^q(B_1)}+ \|u_\l \|_{L^q(B_1)})\leq C ( \|g_\l \|_{L^q(B_1)}+ 1),\nonumber\\
\mbox{ where }q^*=\frac{ qn}{n-2rq}\mbox{ if }2rq <n \mbox{ and for all }   \;q^*> 1 \mbox{ if } q=\frac{n}{2r};
\end{eqnarray}
 and
\begin{eqnarray}\label{3.7'}
\|u_\l\|_{C^{2r-1}(B_\frac 12)} \leq C\|u_\l\|_{W^{2r,q}(B_\frac 12)}\leq C ( \|g_\l \|_{L^q(B_1)}+ \|u_\l \|_{L^q(B_1)})\leq C ( \|g_\l \|_{L^q(B_1)}+ 1),\mbox{ if } 2rq > n.
\end{eqnarray}
So, inequality \eqref{3.7} and  point 4 of Lemma \ref{L11} give
\begin{eqnarray}\label{3.6}
\|g _\l\|_{L^{\frac{q^*}{p_2}}(B_\frac 12)} \leq C \left( \|g_\l \|_{L^q(B_1)}+1\right)^{p_2},\;\mbox{ if }2rq \leq n.
\end{eqnarray}
If $2rq_1\geq n$ (respectively $2rq_1=  n $) the desired  estimate \eqref{zz1} follows from \eqref{ddss} and \eqref{3.7'}
 (with $q=q_1$) (respectively \eqref{ddss}, \eqref{3.7}  (with $q=q_1$) and \eqref{3.7'}, \eqref{3.6} (with $q=p_2\frac{n+1}{2r}$) ).
The case $2rq_1<  n $ needs more involving analysis.  As $q_1=\frac{ p_2+1}{p_2}$ and $1<p_2<\frac{n+2r}{n-2r},$ we have $$q_1^*=\frac{ q_1n}{n-2rq_1}=
\frac{ (p_2+1)n}{p_2(n-2r)-2r}>\frac{ (p_2+1)n}{n}>p_2 \mbox{ and }\frac {1}{q_{1}}-\frac {2rp_2}{n(p_2-1)}<0.$$ Set  $q_2=\frac{q_1^*}{p_2}$ and $q_{k+1}=\frac{q_{k}^*}{p_2}.$ We claim that there exists $k_0\in \N^*$ such that
 \begin{eqnarray}\label{3.8}
 2rq_{k_0+1} >n \mbox{ and } 2rq_{k_0} <n.
 \end{eqnarray}
Suppose by contradiction that $2rq_{k} <n$ for all $k \in \N^*.$ Then, $\frac {1}{q_{k+1}}=\frac {p_2}{q_{k}}-\frac {2rp_2}{n}$ and therefore
\begin{eqnarray}\label{3.11}\frac {1}{q_{k+1}}=\frac {p_2^{k}}{q_{1}}-\frac {2rp_2}{n}\sum_{ j=0}^{k-1}p_2^{j}=p_2^{k}\left(\frac {1}{q_{1}}-\frac {2rp_2}{n(p_2-1)}\right)
  +\frac {2rp_2}{n(p_2-1)}.\end{eqnarray}
  We reach a contradiction since $\frac {1}{q_{k}}\to-\infty.$ Set now
    $$\beta=\frac {2rp_2}{n(p_2-1)}\left( \frac {2rp_2}{n(p_2-1)}-\frac {1}{q_{1}}
    \right)^{-1}=\frac {2r(p_2+1)}{2r(p_2+1)-n(p_2-1)}.$$ From \eqref{3.11}, we have $ p_2^{k_0}< \beta \mbox{  and  } p_2^{k_0+1}>\beta.$
 Hence, iterating \eqref{3.6}, we obtain
   \begin{eqnarray*}
\|g_\l\|_{L^{q_{k_0+1}}(B_{\frac 1{2^{k_0+1}}})} \leq C ( \|g_\l \|_{L^{q_1}(B_1)}+1)^{p_2^{k_0}}\leq C ( \|g_\l \|_{L^{q_1}(B_1)}+1)^{\beta}.
\end{eqnarray*}
Set $\gamma_1=\frac{(p_1+1)\beta}{q_1}=\frac {2r(p_1+1)p_2}{2r(p_2+1)-n(p_2-1)}$ and $\gamma_2= \beta +\frac{2r}{p_1-1}\gamma_1.$ As $rq_{k_0+1}> n,$
the last inequality  with \eqref{3.7'} and \eqref{ddss} imply $$\|u_\l\|_{C^{2r-1}( B_{\frac 1{2^{k_0+1}}})}
\leq C(1+i(u))^{\gamma_2}d_y^{\frac{-2r}{p_1-1}\gamma_1}.$$ According to the definition of $u_\l,$ we get $\sum_{ j=0}^{2r-1}d_y^j|(\nabla^j u)(y)|
\leq  C(1+i(u))^{\gamma_2}d_y^{\frac{-2r}{p_1-1}\gamma_1}.$  This achieves the proof of Theorem \ref{univ}. \qed

 {\bf Appendix A: Proof of Lemma \ref{L11}.} In the following, $C$ denotes generic positive constant depending only on the parameters $(s_0,p_1, p_2)$ and the constant $c_1$
of assumptions $(h_{1})$-$(h_{4})$. The following inequalities are an immediate consequence of $(h_{4})$:
\begin{eqnarray}\label{h4}
|F(x,s)|,\;|f(x,s)s|\leq C, \;\forall \;(x,s)\in \O\times[-s_0,s_0].
\end{eqnarray}
Hence, points 1 and 2 follow from $(h_{1})$-$(h_{2})$. Also, in view of  \eqref{h4} and the fact that the  nonlinearity $-f(x,-s)$ satisfies  $(h_{1})$-$(h_{4}),$  we need only
to prove points 3 and 4 for all $(x,s) \in \O\times [s_0,\infty).$\\
{\bf Proof of point 3.}
  According to $(h_1),$ we have
   \begin{eqnarray}\label{h12}
  f'(x,s)s\geq p_1f(x,s),\,\forall (x,s) \in \O\times [s_0,\infty)
  \end{eqnarray}
 which implies  $\left(\frac{f(x,s) }{s^{p_1}}\right)'\geq 0$. As $f(x,s_0)\geq \frac 1{c_1}$ for all $x\in \O$ (see $(h_4)$), we derive
  \begin{eqnarray}\label{h14}
 f(x,s)\geq \frac{ s^{p_1}}{c_1s_0^{p_1}} \mbox{ and } f(x,s)s\geq \frac{ s^{p_1+1}}{c_1s_0^{p_1}},\,\forall (x,s) \in \O\times [s_0,\infty),
   \end{eqnarray}
which imply the first inequality of point $3$. Integrating now \eqref{h12} over $[s_0,s]$ and using  $(h_2),$ we derive
$\frac{f(x,s)s}{p_2+1}\leq F(x,s)\leq\frac{f(x,s)s}{p_1+1} +C,\,\forall (x,s) \in \O\times [s_0,\infty)$ which pmlies the second and third inequalities  of point $3$.\\
{\bf Proof of point 4.} According to $(h_{2})$, we have $\left(\frac{F(x,s) }{s^{p_2+1}} \right)'\leq 0 \;\forall \;(x,s)\in
  \O\times[s_0,\infty)$  which with \eqref{h4} imply $F(x,s)\ \leq Cs^{p_2+1} \;\forall \;(x,s)\in
  \O\times[s_0,\infty(.$ Hence, from $(h_{2})$ and point 2, we get $|f(x,s)|^{\frac 1{p_2}} \leq C|s|,\;\forall \;(x,s)\in  \O\times[s_0,\infty).$
  Consequently, for $t>0$  and  $q_1=\frac {p_2+1}{ p_2},$ we derive
 \begin{eqnarray*}
 |f(x,s)|^{q_1} \leq C|f(x,s)s| \,\, \mbox{ and }\,\,|f(x,s)|^{\frac {t}{ p_2}}  \leq C|s|^t , \;(x,s)\in \O\times[s_0,\infty).
\end{eqnarray*}

{\bf Proof of point 5.} In view of Young's inequality, we obtain $as^2b\leq \varepsilon s^{p_1+1}a^{\frac{p_1+1}2}+ \varepsilon^{\frac{-2}{p_1-1}}b^{\frac{p_1+1}{p_1-1}}.$ Recall that $0\leq a\leq 1$ and using  point 3, we derive
 point 5. This end the proof of  Lemma \ref{L11}.\qed\\

 {\bf Appendix C: Proof of \eqref{bbcaa}.}\\
Noticing that \eqref{bbcaa} is trivial for $r=1$. Let $k\in\N^*.$ If $r=2k$, i.e. $\textit{D}_{r}=\D^k$, apply Leibnitz's formula, we have
   \begin{eqnarray}\label{bbcaaa}
 \D^k ( \nabla u.\cdot(x-y))=\nabla(\D^k u)\dot(x-y)+2k\D^k u.
 \end{eqnarray}
   Multiplying \eqref{bbcaaa} by $\D^k u$ and taking into account that $\D^k u\nabla(\D^k u)\cdot(x-y)=\frac 12\nabla((\D^k u)^2)\cdot(x-y)$, This achieves the proof of \eqref{bbcaa}.\\
   If $r=2k+1$, that is $\textit{D}_{r}=\nabla\D^k$. According to \eqref{bbcaaa} we derive
 \begin{eqnarray*}
  \textit{D}_{r}u \cdot\textit{D}_{r}(\nabla u\cdot( x-y))= \nabla\D^k\nabla\left(\nabla(\D^k u).( x-y)\right)+(r-1)| \textit{D}_{r}u|^2.
  \end{eqnarray*}
 Therefore,  \eqref{bbcaa} follows as $ \nabla w\cdot\nabla(\nabla w\cdot(x-y))= \frac 12\nabla(|\nabla w|^2)\cdot(x-y)+|\nabla w|^2,\,\forall w\in C^1(\R^n).$ \\

{\bf Appendix C:  Local $L^p$-$W^{2t,p}$-estimate, $t\in \N^*,\,p\geq 2$.} Consider the linear higher order elliptic problem of the form
\begin{equation}\label{EQ}
L u = g\;\; \mbox{ in } \;\; \Omega.
\end{equation}
Here $\O$ is a domain of $\R^n$ and
$$L= \left(- \sum_{i,k=1}^n a_{ik}(x) \frac{\partial^2}{\partial x_i \partial x_k}\right)^r+ \sum_{|j| \leq 2t-1} b_{j} (x) D^{j}$$
 where $ b_{j} \in L_{loc}^{\infty}(\Omega)$, $a_{ik} \in C^{2t-2}(\O)$ and $L$ is a {\bf local} uniformly elliptic operator, that is
For all bounded open subset $O$ of $\O$ there exists a constant $\lambda_{O} >0$ with $\lambda_{O}^{-1} |\xi|^2 \leq \sum_{i,k=1}^n a_{ik}(x) \xi_i \xi_k \leq \lambda_{O} |\xi|^2
$ for all $\xi \in \R^n$, $x \in \overline{O}$. Let $A$
and $A'$ be two bounded open subset of $\O$ such that $\overline{A} \subset A'\subset\overline{A'} \subset \O$ and $\o'.$  When  $p\geq 2$ by virtue of Lemma
\ref{t1}, we establish
  local analogue of the celebrated $L^p$-$W^{2t,p}$ estimate of Agmon-Douglis-Nirenberg \cite{ADN}. Set $d= {\rm dist}(A, \O\backslash A').$
\begin{cor}\label{tloc} Let $g \in L_{loc}^p(\O)$ for some $p  \geq 2.$ Then, there exists a constant $C>0$ depending only on $\|a_{ik} \|_{C^{2r-2}(A')}$, $\|b_{j}\|_{L^\infty(A')}$
 and $\lambda_{A'},A, A',d, n,p,t$ such that for any $u \in W_{loc}^{2t,p}(\O)$ a weak solution of \eqref{EQ}, we have
$$\|u\|_{W^{2t,p}(A)} \leq C \left( \|g \|_{L^p(A')}+ \|u \|_{L^p(A')}\right). $$
\end{cor}

 {\bf Proof of Corollary \ref{tloc}.}\\
 In the following $C$ denotes a generic positive constant which depends on the parameters stated in Corollary \ref{tloc}.  As $\overline{A}$ is a compact subset of \,$A',$ we can
 find $x_i\in A$, $i=1, 2...., i_0$ such that $\overline{A}\subset \bigcup_{1\leq i\leq i_0} B(x_i,\frac{d}{4})\subset A'$ where
$i_0\in \N^*$ { \bf depending only} on $d$ and $\o.$  Let $m\geq 2t$ and $\psi$ be the cut-off function defined in Lemma \ref{00}
 relying on $\o_i=B(x_i,\frac{d}{4})$ and $\o_i'=B(x_i,\frac{d}{2}).$  A simple computations give
\begin{eqnarray*}
L( u\psi^m)= g\psi^m +uL(\psi^m)+b_{0}u\psi^m+\sum_{ 1\leq |i|+|j|\leq 2t-1} c_{i,j} D^juD^i(\psi^m), \mbox{ where } c_{i,j} \in L^\infty(A').
\end{eqnarray*}
As $u\psi^m \in W^{2t,p}(\o')\cap W_0^{t,p}(\o_i')$ with $\o_i'$ is of class $C^{2t}$, Agmon-Douglis-Nirenberg's global
 estimate \cite{ADN}  and \eqref{grad} imply
\begin{eqnarray}\label{loc}
\sum_{ 0\leq s\leq 2t  } \int_{\o'} |\nabla^{s}( u\psi^{m})|^p \leq C\left(\|g\|^p_{L^p(\o')}+\|u\|^p_{L^p(\o')}
+\sum_{1\leq s\leq 2t-1} \sum_{ 1\leq q\leq s  }\int_{\o'} |\nabla^q u|^p|\nabla^{s-q} \psi^m|^p\right).
\end{eqnarray}
 Using now inequality \eqref{ij'} with $r=s$, we obtain
\begin{eqnarray*}
\sum_{ 1\leq q\leq s-1  }\int_{\o'} |\nabla^q u|^p|\nabla^{s-q} \psi^m|^p\leq \varepsilon\int_{\o'} |\nabla^s( u \psi^m)|^p+  C_{\varepsilon, d}\int_{\o'}|u|^p \psi^{p(m-s)}.
\end{eqnarray*}
 Applying again \eqref{ij'} with $r=s+1$ and replacing $\varepsilon$  by $\frac{ \varepsilon}{d}$, yields
 \begin{eqnarray*}
\int_{\o'} |\nabla^s u|^p| \psi^{pm} \leq \int_{\o'} |\nabla^s u|^p \psi^{p(m-1)}\leq \varepsilon \int_{\o'} |\nabla^{s+1}( u \psi^m)|^p+
  C_{\varepsilon, d}\int_{\o'}|u|^p \psi^{p(m-s-1)}.
\end{eqnarray*}
 Collecting the two last inequalities, we derive
 \begin{eqnarray*}
\sum_{1\leq s\leq 2t-1} \sum_{ 1\leq q\leq s  }\int_{\o'} |\nabla^q u|^p|\nabla^{s-q} \psi^m|^p \leq \varepsilon \sum_{ 0\leq s\leq 2t }\int_{\o'} |\nabla^{s}( u \psi^m)|^p+
 C_{\varepsilon, d}\int_{\o'}|u|^p.
\end{eqnarray*}
 We insert the above inequality in the right-hand side of \eqref{loc}and we choose $\varepsilon=\frac 1{2C}$, it follows that
\begin{eqnarray*}
  \| u\psi^{m}\|^p_{W^{2t,p}(\o')} \leq C\left(\|g \|^p_{L^p(\o')}+\|u\|^p_{L^p(\o')} \right).
\end{eqnarray*}
Since $\psi(x)=1$ if $x \in \o,$ we obtain
\begin{eqnarray*}
 \| u\|^p_{W^{2t,p}(B(x_i,\frac{d}{4}))} \leq C\left(\|g \|^p_{L^p(\o')}+\|u\|^p_{L^p(\o')} \right)\leq C\left(\|g \|^p_{L^p(A')}+\|u\|^p_{L^p(A')} \right).
\end{eqnarray*}
As $\overline{A}\subset \bigcup_{1\leq i\leq i_0} B(x_i,\frac{d}{4})$  and $i_0$ depends only on $A$ and $d$, we derive
\begin{eqnarray*}
 \| u\|^p_{W^{2t,p}(A} \leq Ci_0\left(\|g \|^p_{L^p(A')}+\|u\|^p_{L^p(A')} \right).
\end{eqnarray*}
This achieves the proof of Corollary \ref{tloc}.\qed

\end{document}